\documentclass[a4paper,12pt,thmsa]{amsart}
\usepackage[a4paper,marginratio={1:1},scale={0.72,0.74},footskip=7mm,headsep=10mm]{geometry}

\usepackage{amsfonts}
\usepackage{amssymb,amsmath,latexsym}
\usepackage[dvips]{graphics}
\usepackage{graphicx,subfigure}
\usepackage[dvips]{color}
\usepackage[T1]{fontenc}
\usepackage[active]{srcltx}
\usepackage{amsmath}
\usepackage{amsfonts}
\usepackage{amssymb}
\usepackage{psfrag}
\usepackage{color}
\usepackage{url}
\usepackage{amsthm}
\usepackage{array}

\usepackage[english]{babel}
\usepackage{mathrsfs}

\usepackage{pst-tree}
\usepackage{lscape}
\usepackage[T1]{fontenc}
\usepackage{pstricks,pstricks-add}

\newtheorem{thm}{Theorem}[section]
\newtheorem{proposition}[thm]{Proposition}
\newtheorem{remark}[thm]{Remark}
\newtheorem{definition}[thm]{Definition}
\newtheorem{corollary}[thm]{Corollary}

\newtheorem{lemma}[thm]{Lemma}

\newcommand{\e}{\mathbb{E}}

\newcommand{\Pro}{{\bf P}}
\newcommand{\dif}{\mathrm{d}}
\newcommand{\Z}{\mathbb{Z}}

\def\beq{\begin{equation}}               
\def\eeq{\end{equation}}                 
\def\bea{\begin{eqnarray}}             
\def\eea{\end{eqnarray}}               
\def\be*{\begin{eqnarray*}}             
\def\ee*{\end{eqnarray*}}               
\def\ba{\begin{array}}                  
\def\ea{\end{array}}                    

\def\beqlb{\begin{eqnarray}} \def\eeqlb{\end{eqnarray}}
\def\beqnn{\begin{eqnarray*}} \def\eeqnn{\end{eqnarray*}}

\def\<{\langle}  \def\>{\rangle}

\def\bde{\begin{definition}}
\def\ede{\end{definition}}

\def\bth{\begin{theorem}}
\def\eth{\end{theorem}}
\def\bpr{\begin{proposition}}
\def\epr{\end{proposition}}
\def\ble{\begin{lemma}}
\def\ele{\end{lemma}}
\def\bcor{\begin{corollary}}
\def\ecor{\end{corollary}}
\def\bre{\begin{remark}}
\def\ere{\end{remark}}

\def\p{\mathbb{P}}
\def\ee{\varepsilon}

\title[On mutations in the branching model for multitype populations]{On mutations in the branching model\\ 
for multitype populations}

\author{Lo\"ic Chaumont}

\address{Lo\"ic Chaumont -- LAREMA -- UMR CNRS 6093, Universit\'e d'Angers, 2 bd Lavoisier, 49045 Angers cedex~01}

\email{loic.chaumont@univ-angers.fr}

\author{Thi Ngoc Anh Nguyen }

\address{Thi Ngoc Anh Nguyen  -- LAREMA -- UMR CNRS 6093, Universit\'e d'Angers, 2 bd Lavoisier, 49045
Angers cedex~01}

\email{nguyen@math.univ-angers.fr}



\thanks{}

\date{\today}

\begin{document}

\begin{abstract} The forest of mutations associated to a multitype branching forest is obtained by merging 
together all vertices of its clusters and by preserving connections between them. We first show that 
the forest of mutations of any mulitype branching forest is itself a branching forest. Then we give its progeny 
distribution and describe some of its crucial properties in terms the initial progeny distribution. We also obtain 
the limiting behaviour of the number of mutations both when the total number of individuals tends to infinity 
and when the number of roots tends to infinity. The continuous time case is then investigated by considering 
multitype branching forests with edge lengths. When mutations are non reversible, we give a representation 
of their emergence times which allows us to describe the asymptotic behaviour of the latters, when the ratios 
of successive mutation rates tend to 0. 
\end{abstract}

\maketitle

\section{Introduction}

The homogeneous multitype branching hypothesis provides a relevant model of population growth in the 
absence of any competitive or environmental constraint. In particular, it is widely used in population genetics,
when studying successive mutations whose accumulation leads to the development of cancer. Then determining the 
statistics of the emergence times of mutations, or evaluating the distribution of the population size of mutant cells 
at any time become important challenges. In the extensive literature on the subject, let us simply cite \cite{Iwasa2}, 
\cite{Haeno}, \cite{Durrett2010}, and \cite{Durrett}.

This work is concerned with the mathematical study of mutations in multitype branching frameworks. We first 
focus on the problem of the total number of mutations under very general assumptions. This number is not a 
functional of the associated branching process and its study requires the complete knowledge of the multitype 
branching structure, that is the underlying plane forest. Then we show that the forest of mutations associated to any 
multitype forest, is itself a multitype branching forest whose progeny distribution can be 
explicitely computed. This result allows us to investigate the asymptotic behaviour of the number of mutations, when 
either the total population or the initial number of individuals tend to infinity.

When time is continuous, we are mainly interested in emergence times of new mutations in the non reversible 
case. The characterisation of these times requires a good knowledge of the corresponding multitype branching 
process and the main tool in this study consists in a recent extention of the Lamperti representation in higher 
dimensions. Emergence times are then expressed in terms of the underlying multivariate compound Poisson 
process, which allows us to obtain some accurate approximations of their law. 

We start with some preliminaries on the coding of multitype branching forests by multivariate random walks in
Section \ref{9582}. Then we state and prove our results on the total mutations sizes of branching forests in Sections 
\ref{tnm} and \ref{9523}. Results bearing on emergence times are presented in Section \ref{wt}. 
 
\section{Mutations and their asymptotics in discrete multitype forests}\label{5739}

\subsection{Preliminaries on discrete multitype forests}\label{9582}

In all this work, we use the notation $\Z_+=\{0,1,2,\dots\}$ and for any positive integer
$d$, we set $[d]=\{1,\dots,d\}$. We will denote by $e_i$ is the $i$-th unit vector of $\mathbb{Z}_+^d$. We define the
following partial order on $\mathbb{R}^d$ by setting $x=(x_1,\dots,x_d)\ge y=(y_1,\dots,y_d)$, if $x_i\ge y_i$, for all 
$i\in [d]$. The convention $\inf\emptyset=+\infty$ will be valid all along this paper. Then $(\Omega, \mathcal{F}, P)$ is 
a reference probability space on which all the stochastic processes involved in this work are defined.\\

Let us first recall the coding of multitype forests, as it has been defined in \cite{Chau-Liu}. A (plane) forest ${\bf f}$ is a directed 
planar graph with no loops on a possibly infinite and non empty set of vertices ${\bf v} = {\bf v} ({\bf f})$, such that each vertex 
has a finite inner degree and an outer degree equals to 0 or 1. The connected components of a forest are called the {\it trees}. 
In a tree ${\bf t}$, the only vertex with outer degree equal to 0 is called the {\it root} of ${\bf t}$. The roots of the connected 
components of a forest ${\bf f}$ are called the roots of ${\bf f}$. For two vertices 
$u$ and $v$ of a forest ${\bf f}$, if $(u,v)$ is a directed edge of ${\bf f}$, then we say that $u$ is a {\it child} of $v$, or 
that $v$  is the {\it parent} of $u$. We first give an order to the trees of the forest ${\bf f}$ and denote them by 
${\bf t}_1({\bf f}),{\bf t}_2({\bf f}),\dots,{\bf t}_k({\bf f}),\dots$ (we will usually write  ${\bf t}_1,{\bf t}_2,\dots,{\bf t}_k,\dots$ if no 
confusion is possible). Then we rank (a part of) the vertices of ${\bf f}$ according to the breadth first search order, by ranking 
first the vertices of ${\bf t}_1$, then the vertices of ${\bf t}_2$, and so on, see the labeling of the two forests in Figure~\ref{fig0a}. 
Note that if ${\bf t}_k$, for $k\ge1$ is the first infinite tree, then the vertices of ${\bf t}_{k+1},\dots$ have no label according to 
this procedure.\\

To each forest ${\bf f}$, we associate the application $c_{\bf f}:{\bf v}({\bf f})\rightarrow [d]$ such that if 
$u_i,u_{i+1},\dots,$ $u_{i+j}\in{\bf v}({\bf f})$ have the same parent and are placed from left to right, then 
$c_{\bf f}(u_i)\le c_{\bf f}(u_{i+1})\le\dots\le c_{\bf f}(u_{i+j})$. For $v\in{\bf v}({\bf f})$, the integer $c_{\bf f}(v)$ is 
called the {\it type} (or the {\it color}) of $v$. The couple $({\bf f},c_{\bf f})$ is called a {\it $d$-type forest}. When no confusion 
is possible, we will simply write ${\bf f}$. The set of $d$-type forests will be denoted by  $\mathscr{F}_d$.\\

\begin{figure}[hbtp]

\vspace*{-1.3cm}

\begin{center}
{\includegraphics[height=400pt,width=575pt]{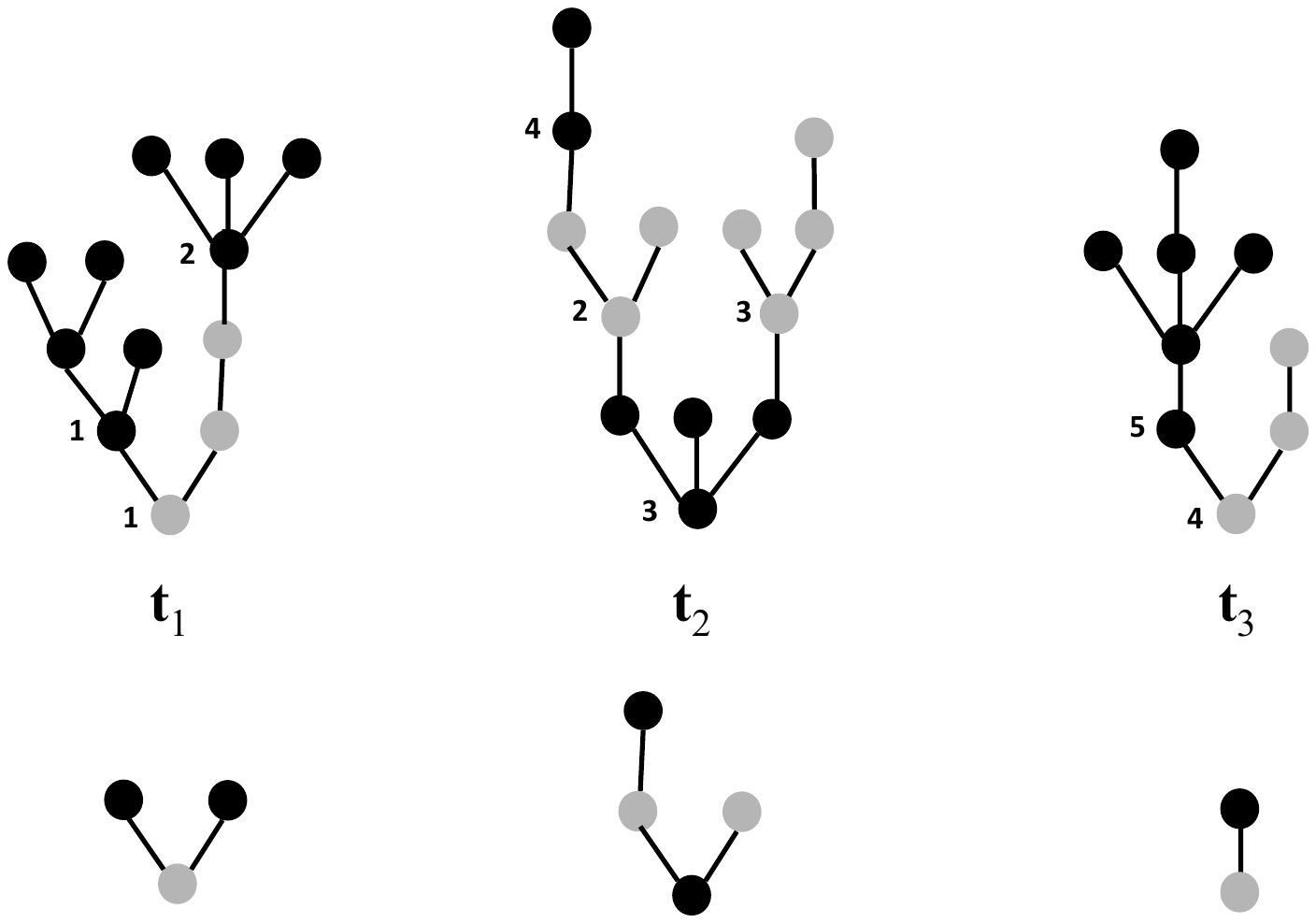}}
\end{center}

\vspace{-6.5cm}

\caption{On the top, a discrete 2-type forest. Roots of clusters are ranked in the breadth first search order of the forest.
The rank is written on the left of these roots. Below, the corresponding forest of mutations.}\label{fig0}
\end{figure}

A {\it cluster} or a {\it subtree of type $i\in[d]$} of a $d$-type forest $({\bf f},c_{\bf f})\in\mathcal{F}_d$ is a maximal connected subgraph 
of $({\bf f},c_{\bf f})$ whose all vertices are of type $i$. Formally, ${\bf t}$ is a cluster of type $i$ of $({\bf f},c_{\bf f})$, if it is a 
connected subgraph whose all vertices are of type $i$ and such that either the root of ${\bf t}$ has no parent or the type of its parent 
is different from $i$. Moreover, if the parent of a vertex $v\in {\bf v}({\bf t})^c$ belongs to ${\bf v}({\bf t})$, then $c_{\bf f}(v)\neq i$. 
Clusters of type $i$ in ${\bf t}_1$ are ranked according to the order of their roots in the breadth first search order of ${\bf t}_1$, see 
Figures \ref{fig0} and \ref{fig0a}. Then if the number of clusters of type $i$ is finite in ${\bf t}_1$, we continue by ranking clusters of 
type $i$ in ${\bf t}_2$, and so on. Note that with this procedure, it is possible that clusters of ${\bf t}_k,{\bf t}_{k+1},\dots$, for some $k$, 
are not ranked. We denote by ${\bf t}^{(i)}_1,{\bf t}^{(i)}_2,\dots,{\bf t}^{(i)}_k,\dots$ the sequence of clusters of type $i$ in 
$({\bf f},c_{\bf f})$. The forest ${\bf f}^{(i)}:=\{{\bf t}^{(i)}_1,{\bf t}^{(i)}_2,\dots,{\bf t}^{(i)}_k,\dots\}$ is 
called {\it the subforest of type $i$} of $({\bf f},c_{\bf f})$. We denote by $u_1^{(i)},u_2^{(i)},\dots$ the elements of ${\bf v}({\bf f}^{(i)})$,
ranked in the breadth first search order of ${\bf f}^{(i)}$. The subforests of the 2-type forest given in Figure~\ref{fig0} are represented 
in Figure~\ref{fig0a}.\\

\begin{figure}[hbtp]

\vspace*{-1cm}

\begin{center}
{\includegraphics[height=400pt,width=575pt]{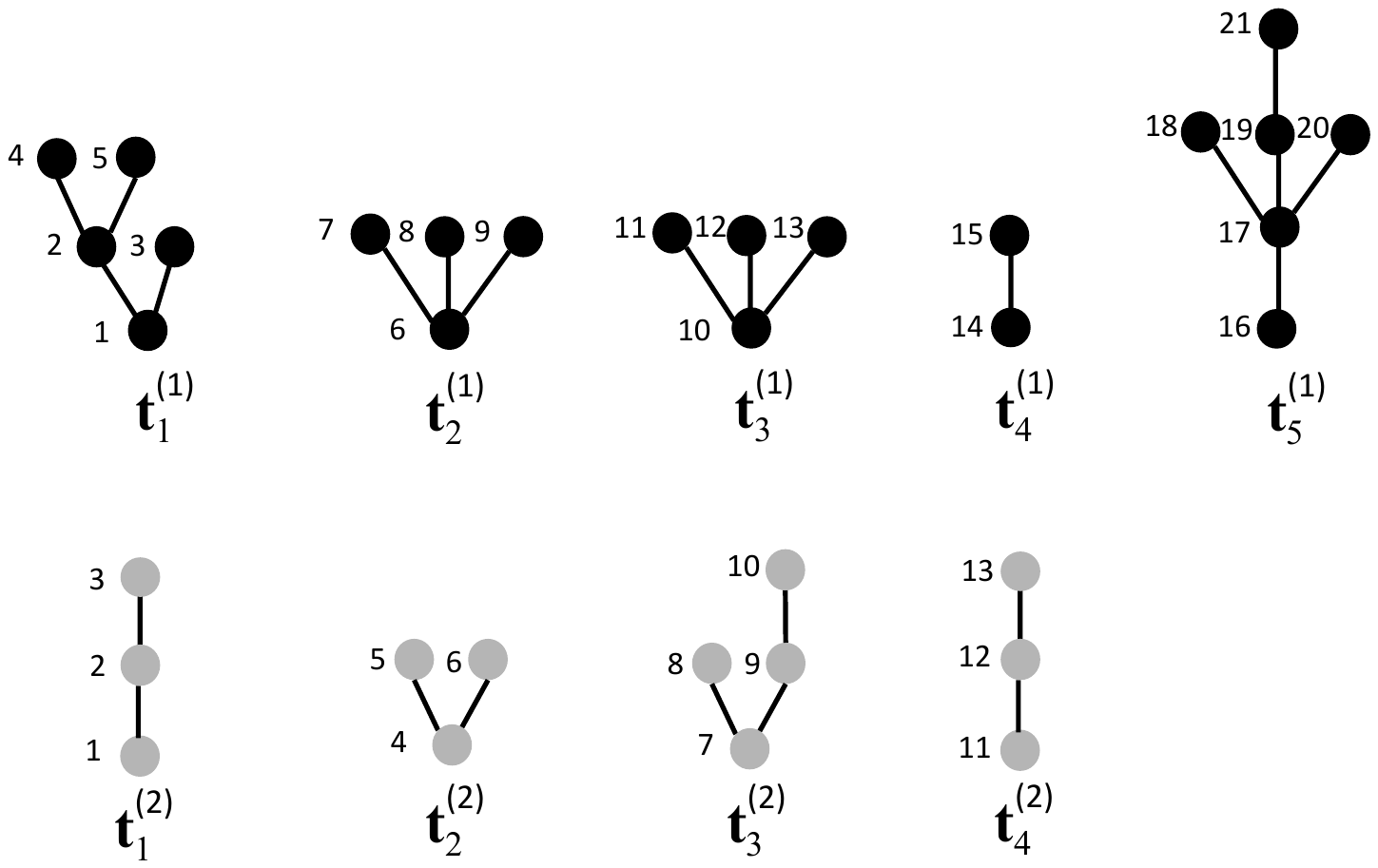}}
\end{center}

\vspace{-6cm}

\caption{The subforests of the 2-type forest given in Figure~\ref{fig0} with their deapth first search labeling.}\label{fig0a}
\end{figure}

To any forest $({\bf f},c_{\bf f})\in\mathcal{F}_d$, we associate the {\it forest of mutations}, denoted by 
$(\bar{{\bf f}},c_{\bar{{\bf f}}})\in\mathcal{F}_d$, which is the forest of $\mathcal{F}_d$ obtained by aggregating all the vertices of each 
subtree of $({\bf f},c_{\bf f})$ with a given type, in a single vertex with the same type, and preserving an edge between each pair 
of connected subtrees. An example is given in Figure~\ref{fig0}.\\

For a forest $({\bf f},c_{\bf f})\in\mathcal{F}_d$ and $u\in{\bf v}({\bf f})$, when no confusion is possible, we denote by $p_i(u)$ the 
number of children of type $i$ of $u$. For each $i\in [d]$, let $n_i\in\mathbb{Z}_+\cup\{\infty\}$ be the number of vertices in the subforest
${\bf f}^{(i)}$ of $({\bf f},c_{\bf f})$. Then let us define the $d$-dimensional chain $x^{(i)}=(x^{i,1},\dots,x^{i,d})$, with length $n_i$ 
and whose values belong to the set $\mathbb{Z}^{d}$, by $x_0^{(i)}=0$ and if $n_i\ge1$,
\begin{equation}\label{6245}
x_{n+1}^{i,j}-x_n^{i,j}=p_j(u_{n+1}^{(i)})\,,\;\;\mbox{if $i\neq j\;$ and}\;\;\;x_{n+1}^{i,i}-x_n^{i,i}=
p_i(u_{n+1}^{(i)})-1\,,\;\;\;0\le n\le n_i-1\,,
\end{equation}
where $(u^{(i)}_n)_{n\ge1}$ is the labeling of the subforest ${\bf f}^{(i)}$ in its own breadth first search order.
Note that the chains $(x_n^{i,j})$, for $i\neq j$ are nondecreasing whereas $(x_n^{i,i})$ is a downward skip free chain, i.e.
$x_{n+1}^{i,i}-x_n^{i,i}\ge -1$, for $0\le n\le n_i-1$. Besides, if $n_i$ is finite, then 
$n_i=\min\{n:x^{i,i}_{n}=\min_{0\le k\le n_i} x^{i,i}_k\}$. Let us also mention that from Theorem 2.7 of \cite{Chau-Liu}, when trees 
of $({\bf f},c_{\bf f})$ are finite, the data of the chains $x^{(1)},\dots,x^{(d)}$ together with the sequence of ranked roots of 
$({\bf f},c_{\bf f})$, allow us to reconstruct this forest.\\

Let us now apply this coding to multitype branching forests. Let $\nu:=(\nu_1,\dots,\nu_d)$, where $\nu_i$ is some distribution 
on $\mathbb{Z}_+^d$. We consider a branching process with progeny distribution $\nu$, that is a population of individuals 
which reproduce independently of each other at each generation. Individuals of type $i$ give birth to $n_j$ children of type $j\in[d]$ 
with probability $\nu_i(n_1,\dots,n_d)$.  For $i,j\in[d]$, we denote by $m_{ij}$ the mean number of children of type $j$, given by an 
individual of type $i$, i.e.
\[m_{ij}=\sum_{(n_1,\dots,n_d)\in\mathbb Z_+^d}n_j\nu_i(n_1,\dots,n_d)\,.\]
We say that $\nu$ is non singular if there is $i\in[d]$ such that $\nu_i({\rm n}:n_1+\dots+n_d=1)<1$. The matrix $M=(m_{ij})$ is said 
to be irreducible if for all $i,j$, $m_{ij}<\infty$ and there exists $n\ge1$ such that $m_{ij}^{(n)}>0$, where $m_{ij}^{(n)}$ is the $ij$ 
entry of the matrix $M^n$. If moreover the power $n$ does not depend on $(i,j)$, then $M$ is said to be primitive. In the latter case, 
according to Perron-Frobenius theory, the spectral radius $\rho$ of $M$ is the unique eigenvalue which is positive, simple and 
with maximal modulus. If $\rho\le1$, then the population will become extinct almost surely, whereas if $\rho>1$, then with positive
probability, the population will never become extinct. We say that $\nu$ is subcritical if $\rho<1$, critical if $\rho=1$ and supercritical 
if $\rho>1$. We sometimes say that $\mu$ is irreducible, primitive, (sub)critical or supercritical, when this is the case for $M$.\\

By multitype branching forest with progeny distribution $\nu$, we mean a sequence with a finite (deterministic) or infinite number 
of independent multitype branching trees with progeny distribution $\nu$. A multitype branching forest will be considered as a random 
variable defined on the probability space $(\Omega, \mathcal{F}, P)$ and with values in $\mathscr{F}_d$. To any multitype branching 
forest $F$, we associate the random sequences $X=\{X^{(i)},i\in[d]\}$, where $X^{(i)}=\{(X^{i,1}_n,\dots,X^{i,d}_n),0\le n\le n_i\}$, which 
are constructed as in (\ref{6245}). It has been proved in \cite{Chau-Liu}, Theorem 3.1 that if $F$ is a primitive and (sub)critical 
branching forest with a finite number of trees, then $X^{(i)}$, $i\in[d]$ are independent random walks whose step distribution 
$\tilde{\nu}_i$ is defined by
\begin{equation}\label{4590}
\tilde{\nu}_i(k_1,\dots,k_d):=\nu_i(k_{1},\dots,k_{i-1},k_{i}+1,k_{i+1},\dots,k_{d})\,,\;\;\mbox{for all 
$(k_1,\dots,k_d)\in\mathbb{Z}_+^d$,}
\end{equation}
and stopped at the smallest solution $(N_1,\dots,N_d)$ of the system
\begin{equation}\label{7321}
x_j+\sum_{i=1}^d X^{i,j}(N_i)=0\,,\;\;\;j\in[d]\,.
\end{equation}
In this equation, $N_i$ is the total number of vertices of type $i$ in $F$ and $x_i$ is the total number of trees in this forest whose root 
is of type $i$. We will say that $F$ is issued from $x=(x_1,\dots,x_d)$. Note that the variables $N_i$ are random, whereas the $x_i$'s 
are deterministic.

\subsection{The total number of mutations and its asymptotics}\label{tnm}

A mutation of type $i$, is the birth event of an individual of type $i$ from an individual of any type $j\neq i$. 
The aim of this section is to study the evolution of mutations in a multitype branching forest. 
Our main result asserts that the forest of mutations, that is the forest obtained by merging 
together all the vertices of a same cluster, is itself a branching forest if and only if for each $i\in[d]$, one of the 
following conditions is satisfied,
\begin{eqnarray*}
&&(A_i)\qquad m_{ii}\le 1\,,\\
&&(B_i)\qquad \mbox{$m_{ii}>1$ and for all $j\neq i$, $m_{ij}=0$.}
\end{eqnarray*}
Moreover, its progeny distribution can be expressed in terms of this of the initial forest. Note that the branching property of 
the forest of mutations is intuitively clear. In the neutral case, it has been pointed out in \cite{Taib}. 

\begin{thm}\label{2358} Let $F$ be any multitype branching forest with progeny distribution $\nu=(\nu_1,\dots,\nu_d)$ and denote by 
$\overline{F}$ the associated forest of mutations. Assume that for all $i\in[d]$, one of the conditions $(A_i)$ or $(B_i)$ holds. 
Then $\overline{F}$ is a multitype branching forest with progeny distribution $\mu=(\mu_1,\dots,\mu_d)$ on $\mathcal{S}_i:=\{{\rm k}\in\mathbb{Z}_+^d:k_i=0\}$, which is defined by
\begin{equation}\label{4146}
\mu_i({\rm k})=\sum_{n\ge1}n^{-1}\nu^{*n}_i({\rm k}+(n-1)e_i)\,,\;\;\;{\rm k}\in\mathcal{S}_i\,,
\end{equation}
if $(A_i)$ is satisfied. If $(B_i)$ is satisfied, then $\mu_i$ is the Dirac mass at $0$. Moreover $\mu$ satisfies the following properties:
\begin{itemize}
\item[$1.$] Let $\overline{M}=(\bar{m}_{ij})$ be the mean matrix of $\mu$ and let $r\ge1$. Then $\mu_i$ admits moments of order $r$ 
if and only if either for all $j\neq i$, $m_{ij}=0$ or $\nu_i$ admits moments of order $r$ and $m_{ii}<1$. In the latter case, for all $i,j$ 
such that $i\neq j$, $\bar{m}_{ij}=\frac{m_{ij}}{1-m_{ii}}$.
\item[$2.$]  Assume that $\bar{m}_{ij}<\infty$, for all $i,j\in[d]$. Then $\overline{M}$ is irreducible if and only if $M$ is irreducible. 
If $\overline{M}$ is primitive, then so is $M$. The converse is not true. 
\item[$3.$] Assume that $\overline{M}$ is primitive, then $\overline{M}$ is subcritical $($resp. critical, resp. super-critical$)$ if and 
only if $M$ is subcritical $($resp. critical, resp. supercritical$)$.
\end{itemize}
If for some $i\in[d]$, none of the conditions $(A_i)$ and $(B_i)$ holds, then there is $j\neq i$ such that individuals of type $i$ 
in $\overline{F}$ give birth to an infinite number of children of type $j$ with positive probability. Therefore $\overline{F}$ is not a 
branching forest in our sense. 
\end{thm}
\begin{proof}  Since the result only bears on the progeny law of forests, we do not loose any generality by assuming that $F$ has an 
infinite number of trees. Then the stochastic processes $X=\{X^{(i)},\,i\in[d]\}$ obtained from $F$, as in (\ref{6245}) are defined on
the whole integer line $\{0,1,\dots\}$. 
Note that their definition slighly extends the definition which is given in \cite{Chau-Liu}. Indeed, without any more assumption on $\nu$, 
trees of the forest can be infinite, so that the process $X$ is not necessarily a coding of the forest, that is, if some trees are infinite then 
it is not  possible to reconstruct the whole forest from $X$ and the sequence of its roots. However, it is straightforward to check that 
$X^{(i)}$, $i\in[d]$ are independent random walks and that the step distribution of $X^{(i)}$ is $\tilde{\nu}_i$, which is defined in 
(\ref{4590}). In particular, the law of $X$ characterizes this of $F$.  

Now, let us consider the forest of mutations $\overline{F}$. By construction, this forest is composed of an infinite number of 
independent and identically distributed trees. Hence, in order to show that $\overline{F}$ is a branching forest, it suffices to 
show that its trees are branching trees. 

Let us denote by $\{\overline{X}^{(i)},i\in[d]\}$ the process which is defined from $\overline{F}$ as in (\ref{6245}). Let $i\in[d]$ and 
assume first that $(A_i)$ holds. Then we define the first passage time process of the random walks $X^{i,i}$, $i\in[d]$ by,
\[\tau_k^{(i)}=\inf\{n\ge0:X^{i,i}_n=-k\}\,,\;\;\;k\ge0\,.\]
Since $m_{ii}\le1$, then from the law of large numbers, $\liminf_{n\rightarrow\infty}X^{i,i}_n=-\infty$, a.s., so that $\tau_k^{(i)}$ is 
almost surely finite for all $k\ge0$ and $\lim_{k\rightarrow\infty}\tau^{(i)}_k=\infty$, a.s. Moreover, for all $i,j\in[d]$, 
\[\overline{X}^{i,j}_k=X^{i,j}(\tau_k^{(i)})\,,\;\;\;k\ge0\,.\]
Indeed, the effect of the time change by $\tau_k^{(i)}$ is to merge all vertices of a same cluster of type $i$ into a single vertex. 
Note that $\overline{X}^{(i)}$, $i\in[d]$ are independent random walks. Assume with no loss of generallity that the root of the first 
tree in $\overline{F}$ has type $1$, then a slight extention Theorems 2.7 and 3.1 in \cite{Chau-Liu} to any progeny distribution, 
allows us to show that this first tree is coded by the processes $(\overline{X}^{(i)}_k,0\le k\le N_i)$, $i\in[d]$, where $(N_1,\dots,N_d)$ 
is the smallest solution of the system
\begin{equation}\label{3721}
r_j+\sum_{i=1}^d \overline{X}^{i,j}(N_i)=0\,,\;\;\;j\in[d]\,,
\end{equation}
and $(r_1,\dots,r_d)=(1,0,\dots,0)$. Note that in our case, $N_i$ can be infinite. This extended notion of smallest solution is defined 
in \cite{Chaumont}, see Lemma 1 therein. This coding result implies that the first tree in $\overline{F}$ can be reconstructed from the 
processes $(\overline{X}^{(i)}_k,0\le k\le N_i)$, $i\in[d]$ and applying part 3. of Theorem 3.1 in \cite{Chau-Liu}, we obtain that this 
tree is a branching tree whose progeny distribution $\mu=(\mu_i,i\in[d])$ is given by
\[\mu_i(k_1,\dots,k_d)=P(\overline{X}^{(i)}_1=(k_1,\dots,k_{i-1},-1,k_{i+1},k_d))\,,\;\;\;(k_1,\dots,k_d)\in\mathcal{S}_i\,.\]
Then in order to make this law explicit in terms of $\nu$, we apply the Ballot theorem for cyclically exangeable sequences due to 
Tak\'acs \cite{ta}. Since conditionally on $X^{i,j}$, $i\neq j$, $X^{i,i}$ is downward skip free with cyclical exchangeable increments, we 
have for all $(k_1,\dots,k_d)\in\mathcal{S}_i$,
\begin{eqnarray*}
&&P(\overline{X}^{(i)}_1=(k_1,\dots,k_{i-1},-1,k_{i+1},\dots,k_d))\\
&=&\sum_{n\ge1}
P(X^{(i)}_n=(k_1,\dots,k_{i-1},-1,k_{i+1},\dots,k_d),\tau_1^{(i)}=n)\\
&=&\sum_{n\ge1} \frac1nP(X^{(i)}_n=(k_1,\dots,k_{i-1},-1,k_{i+1},\dots,k_d))\,,
\end{eqnarray*}
which gives (\ref{4146}) from (\ref{4590}).
If $(B_i)$ holds, then by definition, individuals of type $i$ in $\overline{F}$ are all leaves and hence,
$\overline{X}^{i,j}\equiv0$, for all $j\neq i$ and $\overline{X}^{i,i}_n=-n$, for all $n\ge0$, see (\ref{6245}). In this case, the conclusion
follows immediately. 

Let us now prove properties 1--3 of $\mu$. First note that for all $i\neq j$, $m_{ij}=0$ if and only if $\bar{m}_{ij}=0$. Then let $r\ge1$,
assume that $\mu_i$ admits moments of order $r$ and that there is $j\neq i$ such that $m_{ij}=E(X^{i,j}_1)>0$. The variable 
$\tau_1^{(i)}$ is a stopping time in the filtration generated 
by $X^{(i)}$ to which the increasing random walk $X^{i,j}$ is adapted. Then by applying Theorem 5.4 in \cite{gu}, we obtain that
$E((X^{i,j}_1)^r)<\infty$ and $E((\tau_1^{(i)})^r)<\infty$. In particular $\tau_1^{(i)}<\infty$, a.s. Now by definition, the random walk
$(X^{i,i}_n)$ can be written as $X^{i,i}_n=Y^{i,i}_n-n$, where $(Y_n^{i,i})$ is an increasing random walk. Since 
$Y^{i,i}(\tau_1^{(i)})=\tau_1^{(i)}-1$ and $E((\tau_1^{(i)})^r)<\infty$, we have $E\left(|Y^{i,i}(\tau_1^{(i)})|^r\right)<\infty$ and by 
applying Theorem 5.4 in \cite{gu} again, we obtain that $E\left(|Y^{i,i}_1|^r\right)<\infty$, and hence $E\left(|X^{i,i}_1|^r\right)<\infty$. 
So we have proved that $\nu$ admits moments of order $r$. Then it follows from the definition of $\tau^{(i)}_1$ and from Lemma 3.1 
in \cite{km} that $E((\tau^{(i)}_1)^r)<\infty$ implies $\lim_{n\rightarrow\infty} X^{i,i}_n=-\infty$, and hence $m_{ii}<1$, from the law of 
large numbers.

Conversely, if  $m_{ij}=0$ for all $j\neq i$, then $\bar{m}_{ij}=0$ for all $j\neq i$ and $\mu_i$ is the Dirac mass at 0, so it 
admits moments of order $r$. Now assume that $\nu_i$ admits moments of 
order $r$ and $m_{ii}<1$. Then it follows directly from Lemma 3.1 in \cite{km} that $E((\tau^{(i)}_1)^r)<\infty$. Moreover from 
Theorem 5.2 in \cite{gu}, 
$E(X^{i,j}(\tau_1^{(i)})^r)<\infty$, for all $j\neq i$, which means that $\mu_i$ admits moments of order $r$. If $\nu_i$ admits 
moments of order 1 and $m_{ii}<1$, then it follows from the optional stopping theorem applied to the martingale 
$(X^{i,j}_n-nE(X^{i,j}_1))$, that $E(X^{i,i}(\tau^{(i)}_1))=-1=E(X^{i,i}_1)E(\tau_1^{(i)})=(m_{ii}-1)E(\tau_1^{(i)})$, and when 
$i\neq j$, $E(X^{i,j}(\tau^{(i)}_1))=E(X^{i,j}_1)E(\tau_1^{(i)})=\frac{m_{ij}}{1-m_{ii}}$ and part 1 is proved. 

If $\overline{M}$ is irreducible, then for all $i$, there is $j\neq i$ such that $\bar{m}_{ij}>0$. From part 1., $\nu_i$ admits moments 
of order 1 and $m_{ii}<1$, for all $i$. In this case, 
\[\overline{M}+\Delta_2=\Delta_1M\,,\;\;\mbox{where}\;\;\Delta_1=\mbox{diag}(\frac{1}{1-m_{ii}})\;\;\mbox{and}\;\;
\Delta_2=\mbox{diag}(\frac{m_{ii}}{1-m_{ii}})\,,\]
and we derive from this identity that $M$ is irreducible. Conversely if $M$ is irreducible, then for all $i$, there is $j\neq i$ such that
$m_{ij}>0$ and hence $\bar{m}_{ij}>0$. Since by assumption, $\bar{m}_{ij}<\infty$, for all $i,j$, then from part 1., $m_{ii}<1$, and
$\overline{M}+\Delta_2=\Delta_1M$ holds. We derive from this identity that $\overline{M}$ is irreducible. 

Now if $\overline{M}$ is primitive, then it is irreducible and as before, $m_{ii}<1$ for all $i\in[d]$. Moreover, 
\[M=(I-\mbox{diag}(m_{ii}))\overline{M}+\mbox{diag}(m_{ii})\,.\] 
Therefore $M$ is primitive. The converse cannot be true since there are nonnegative, irreducible matrices whose main diagonal is zero
and which are not primitive. We can find distributions $\nu$ such that it is the case for $\overline{M}$ and hence for $(I-\mbox{diag}(m_{ii}))\overline{M}$. If $m_{ii}>0$, for all $i$, then it follows from general theory of nonnegative matrices that
 $M=(I-\mbox{diag}(m_{ii}))\overline{M}+\mbox{diag}(m_{ii})$ becomes primitive, see \cite{se}. 

Let us now prove 3. Recall that by definition, since $\overline{M}$ is primitive, $\mu_i$ admits moments of order 1 for all $i\in[d]$. 
Then from the same arguments as in part 2., $M=(I-\mbox{diag}(m_{ii}))\overline{M}+\mbox{diag}(m_{ii})$ 
and $m_{ii}<1$ for all $i\in[d]$. Assume that $M$ is surpercritical, then there is a positive vector $x$ such that $Mx>x$. Therefore, 
$(I-\mbox{diag}(m_{ii}))\overline{M}x>(I-\mbox{diag}(m_{ii}))x$ and since $m_{ii}<1$, we obtain $\overline{M}x>x$. Hence 
$\overline{M}$ is supercritical. Conversely, assume that $\overline{M}$ is supercritical. Then there is a positive vector $x$ such that 
$\overline{M}x>x$, so that $Mx=(I-\mbox{diag}(m_{ii}))\overline{M}x+\mbox{diag}(m_{ii})x>(I-\mbox{diag}(m_{ii}))x+
\mbox{diag}(m_{ii})x=x$ and thus $M$ is supercritical.  Then the identity $M=(I-\mbox{diag}(m_{ii}))\overline{M}+
\mbox{diag}(m_{ii})$ allows us to derive that $M$ is critical if and only if this is the case for $\overline{M}$.

Finally assume that $m_{ii}>1$ for some $i\in[d]$. If $m_{ij}=0$, for all $j\neq i$, then it is clear that individuals of type $i$ in 
$\overline{F}$ are leaves. If $m_{ij}>0$, for some $j\in[d]$, then since clusters of type $i$ are supercritical, some of them have 
infinitely many children with positive probability. Conditionally to this event, such a cluster produces almost surely infinitely many children 
of type $j$, which is equivalent to say that individuals of type $i$ in $\overline{F}$ give birth to an infinite number of children of type $j$ 
with positive probability. 
$\;\;\;\Box$\\
\end{proof}

Let us now consider a multitype branching forest $F$ with progeny distribution $\nu$, with a finite number of trees and let 
$Z_n=(Z^{(1)}_n,\dots,Z^{(d)}_n)$, $n\ge0$ be the associated branching process, that is for each $i\in[d]$, $Z^{(i)}_n$ is the total 
number of individuals of type $i$ present in $F$ at generation $n$.  For $x=(x_1,\dots,x_d)\in\mathbb{Z}_+^d$, we denote by $\Pro_x$ 
the law on $(\Omega,\mathcal{F})$ under which $F$ is issued from $x$. In particular, $\Pro_x(Z_0=x)=1$. Then the next result gives 
the law of the total number of mutations in the forest $F$, that is the number of mutations up to the last generation whose rank is
the extinction time, $T:=\inf\{n:Z_n=0\}$. For $i,j\in[d]$, let us denote by $M_i$  the total number of mutations of type $i$ in $F$, up to 
time $T$ and by $M_{ij}$ the total number of mutations of type $j$ produced by individuals of type $i$. In particular, $M_{ii}=0$ and 
$M_i$ and $M_{ij}$ satisfy the relations
\[M_j=\sum_{i=1}^d M_{ij}\,,\;\;\;j\in[d].\]
Note that if $\nu$ is primitive and supercritical, then $\Pro_x(T=\infty)>0$ for all $x\in\mathbb{Z}_+^d$, so that under $\Pro_x$, $M_i$
and $M_{ij}$ are infinite with positive probability, for some $i,j\in[d]$. We also emphasize that $M_i$ and $M_{ij}$ are not functionals of 
the branching process $(Z_n)$.

\begin{corollary}\label{7279} Assume that $(A_i)$ or $(B_i)$ holds for all $i\in[d]$. Then for all integers $x_i, n_i, k_{ij}, i,j\in [d]$, such 
that $x_i\ge0$, $n_i=-k_{ii}$, for $i\neq j$, $k_{ij}\geq 0$, and for all $j\in[d]$, $n_j=x_j+\sum_{i\neq j}k_{ij} $,
\begin{align*}
&\Pro_x\left(M_1=n_1-x_1,\dots,M_d=n_d-x_d, M_{ij}=k_{ij}, \forall i\neq j\right)\\
&\quad=\frac{\det(K)}{\bar{n}_1\dots\bar{n}_d}
\prod_{i=1}^d\mu_i^{* n_i}(k_{i1},\dots,k_{i(i-1)},0,k_{i(i+1),\dots,k_{id}}),
\end{align*}
where $\mu_i$ is defined in Theorem $\ref{2358}$ and
$\mu_i^{* 0}=\delta_0, \bar{n}_i=n_i\vee 1$,  $K$ is the matrix $(-k_{ij})_{i,j}$ to which we removed 
the line $i$ and the column $i$ for all $i$ such that $n_i=0$.
\end{corollary}
\begin{proof}
This result is a direct consequence of Theorem 1.2 in \cite{Chau-Liu} and Theorem \ref{2358} applied to the forest of mutations
$\overline{F}$. Indeed, it suffices to note that $x_i+M_i$ corresponds to the total number of individuals of type $i$ in $\overline{F}$.
Note however that Theorem 1.2 in \cite{Chau-Liu} is proved only in the case where $\nu$ is primitive and (sub)critical.
But using the coding which is presented in Section \ref{9582} and appyling Lemma 1 in \cite{Chaumont}, we can 
check that it is still valid in the general case by following along the lines the proof which is given in \cite{Chau-Liu}.
$\;\;\;\Box$\\
\end{proof}
\noindent If for some $i\in[d]$, none of the conditions $(A_i)$ and $(B_i)$ holds, then the definition of the vector of mutation sizes 
$(M_1,\dots,M_d)$ still makes sense. In this case, it is possible to obtain its law by extending Theorem \ref{2358} to branching 
forests whose progeny laws give mass to infinity. Note also that Corollary \ref{7279} can be considered as an extension of 
Theorem 1 in \cite{Bertoin}, where a similar formula is given in the neutral case.\\

We now turn our attention to the asymptotic behaviour of the number of mutations, when the total population is growing to infinity.
Our first result is concerned with the critical case and is a direct consequence of Proposition 2 in \cite{pe} and Theorem \ref{2358}. 
If $M$ is primitive, then we denote by $u$ and $v$ the unique right and left positive eigenvectors of $M$ which are associated to 
the eigenvalue 1 and normalized by $u.1=u.v=1$. Recall that, for a multitype branching forest $F$, when no confusion is possible, 
$N_i$ denotes the total population of type $i$ in $F$ and $M_i$ denotes the total number of mutations of type $i$ in $F$. Note also
that when $\nu$ is primitive and critical, then $(A_i)$ necessarily holds for all $i\in[d]$, so that from Theorem \ref{2358}, the forest 
of mutations $\overline{F}$ associated to $F$ is a branching forest with progeny distribution $\mu$ defined by (\ref{4146}). 

\begin{corollary} Let $F$ be a branching forest with a non singular, primitive and critical progeny distribution $\nu$. Assume that for 
all $i\in[d]$, $\mu_i$ admits moments of order $d+1$.  If moreover $\overline{M}$ is primitive and 
the covariance matrices $\Sigma^i$, $\overline{\Sigma}^i$ of $\nu_i$ and $\mu_i$, respectively are positive definite. Then $m_{ii}<1$, 
for all $i\in[d]$ and there are constants $C_1,C_2>0$ such that for all $x_0\in\mathbb{Z}_+^d$,
\begin{eqnarray*}
\lim_{n\rightarrow\infty}n^{d/2+1}\Pro_{x_0}(M_i=\lfloor n(1-m_{ii})v_i\rfloor,\; i\in[d])&=&C_1x_0.u\,,\\\label{5728}
\lim_{n\rightarrow\infty}n^{d+1}\Pro_{x_0}(M_i=\lfloor n(1-m_{ii})v_i\rfloor,\;N_i=\lfloor nv_i\rfloor\,,\; i\in[d])&=&
C_2x_0.u\,.\label{5728}
\end{eqnarray*}
\end{corollary}
\begin{proof}  Since by assumption, $\overline{M}$ is primitive, then for all $i$, there is $j\neq i$ such that $\bar{m}_{ij}>0$, and 
hence $m_{ij}>0$. Therefore, from part $1.$ of Theorem \ref{2358}, $m_{ii}<1$, for all $i$. 
Moreover, from our assumptions and part 3. of Theorem \ref{2358}, $\mu$ is critical. Besides, it is plain that $\overline{M}$ is non 
singular. Then conditions of Proposition 2 in \cite{pe} are satisfied for the multitype branching process associated to $\overline{F}$ 
and the first assertion follows with $\bar{u}$ and $\bar{v}$, the normalized, positive right and left eigenvectors of $\overline{M}$ 
associated to the eigenvalue 1.  Then recall from the proof of part 3. of Theorem \ref{2358} that 
$M=(I-\mbox{diag}(m_{ii}))\overline{M}+\mbox{diag}(m_{ii})$. We derive from this identity that $\bar{u}=u$ and 
$\bar{v}=cv(I-\mbox{diag}(m_{ii}))$, where $c=\|u\cdot v(I-\mbox{diag}(m_{ii}))\|^{-1}$ and the first assertion follows. 

The proof of the second assertion follows the same lines as the proof of Proposition 2 in \cite{pe}. In this case, since the number
of mutations is taken into account together with the total number of individuals, a $2d$-dimensional random walk is involved in 
the proof, which explains that the rate of convergence in now $d+1$. 
$\;\;\;\Box$\\
\end{proof}
\noindent Note that the constants $C_1$ and $C_2$ can be made explicit in terms of the distributions $\nu$ and $\mu$ by 
properly exploiting the proof of Proposition 2 in \cite{pe}.\\

Through the next result we focus on the asymptotic behaviour of the number of mutations in a branching forest when the 
initial number of individuals $x=(x_1,\dots,x_d)$ tends to infinity along some given direction.  

\begin{thm} Let $F(x)$ be any family of multitype branching forests defined on the space $(\Omega,\mathcal{F},P)$, indexed by
$x\in\mathbb{Z}_+^d$ and such that for each $x$, $F(x)$ has progeny distribution $\nu$ and is issued from $x$. 
For $i\in[d]$, let $N_i(x)$ $($resp. $M_i(x)$$)$ be the total number of individuals $($resp. of mutations$)$ of type $i$ 
in $F(x)$. Assume that $\nu$ is primitive and let $w\in\mathbb{Z}_+^d\setminus \{0\}$. 
\begin{itemize}
\item[$1.$] If $\nu$ is critical, then 
\[\lim_{n\rightarrow\infty}\frac{N_i(nw)}{n}=\infty\;\;\mbox{and}\;\;
\lim_{n\rightarrow\infty}\frac{M_i(nw)}{N_i(nw)}=1-m_{ii}\,,\;\;\mbox{in probability.}\]
\item[$2.$] If $\nu$ is subcritical, then 
\[\lim_{n\rightarrow\infty}\frac{N_i(nw)}{n}=c_i(w)\;\;\mbox{and}\;\;
\lim_{n\rightarrow\infty}\frac{M_i(nw)}{n}=w_i+(1-m_{ii})c_i(w)\,,\;\;\mbox{in probability,}\]
where $c_i(w):=\sum_{k=1}^dw_k(I-M)^{-1}_{ki}$.
\end{itemize}
In any case, $m_{ii}<1$, for all $i\in[d]$.
\end{thm}
\begin{proof} In order to prove our result, it suffices to construct some particular family of forests $F(x)$, such that for each $x$, 
$F(x)$ has progeny distribution $\nu$ and is issued from $x\in\mathbb{Z}_+^d$, and to show that the limits in the statement hold. 
 
Recall the coding of multitype branching forests which is presented at the end of Section \ref{tnm} and let $X^{(i)}=\{X^{i,j},j\in[d]\}$
be $d$ independent random walks whose respective step distributions are $\tilde{\nu}_i$, $i\in[d]$ defined in (\ref{4590}). Then for 
each $x\in\mathbb{Z}_+^d$, we construct a forest $F(x)$ such that $F(x)$ is encoded by the random walks $X^{(i)}$, $i\in[d]$ and 
contains exactly $x_i$ trees whose root is of type $i$. This construction is possible in the primitive, (sub)critical case, thanks to part 
3. of Theorem 3.1 in \cite{Chau-Liu}.

Then $N_i(x)$ and $X^{(i)}$, $i\in[d]$, satisfy identity (\ref{7321}). Moreover, for $k\neq i$, the number of mutations of type $i$ issued 
from an individual of type $k$ is $X^{k,i}(N_k(x))$, so that the total number of mutations of type $i$ is 
\[M_i(x)=\sum_{k\neq i} X^{k,i}(N_k(x))=-x_i-X^{i,i}(N_i(x))\,.\]
We derive from Lemma 2.2 in \cite{Chau-Liu}, that if $x_1,x_2\in\mathbb{Z}_+^d$ are such that $x_1\le x_2$, 
then the couple of random variables $(N_i(x_2)-N_i(x_1),X^{i,i}(N_i(x_2)-X^{i,i}(N_i(x_1)))$ is independent of process 
$((N_i(x),X^{i,i}(N_i(x))),x\le x_1)$ and has the same law as 
$(N_i(x_2-x_1),X^{i,i}(N_i(x_2-x_1))$. Therefore, for any $w\in\mathbb{Z}_+^d$, $((N_i(nw),X^{i,i}(N_i(nw)),n\ge0)$ is a bivariate
random walk whose step distribution is the law of $(N_i(w),X^{i,i}(N_i(w))$. 

Let $Z=(Z^{(1)},\dots,Z^{(d)})$ be the branching process associated to $F(w)$. Then by definition of $N_i(w)$,  
we have $N_i(w)=\sum_{n=0}^\infty Z^{(i)}_n$.  But ${\bf E}_{w}(Z_n)=wM^n$, so that 
${\bf E}_{w}(Z^{(j)}_n)=\sum_{i=1}^dw_im_{ij}^{(n)}$ and since $\nu$ is primitive, we have from Frobenius Theorem for primitive 
matrices, $m_{ij}^{(n)}\sim u_iv_j\rho^n$, see Theorem 1, Section V.2 in \cite{Athreya}. 
So we have proved that $E(N_i(w))<\infty$ if and only if $\nu$ is subcritical. Moreover, if $\nu$ is subcritical, then $I-M$ is invertible 
and it follows from the above expressions that $E(N_i(w))=\sum_{i=1}^dw_i(I-M)^{-1}_{ij}$. Then assertions 1. and 2. follow directly 
from the law of large numbers.  

Finally, since $\nu$ and $\mu$ are primitive, by definition, they admit moments of order 1 and we derive from part 1. of Theorem 
\ref{2358} that $m_{ii}<1$, for all $i\in[d]$.$\;\;\;\Box$\\
\end{proof}

\section{When continuous time is involved}\label{2682}

\subsection{The Lamperti representation}\label{3249}

Let us now consider a $d$ type population which is composed at time $t=0$, of $x_i$ individuals of 
type $i\in[d]$ and whose dynamics in continuous time behave according to a branching model.
More specifically, at any time, all individuals in the population live, give birth and die independently of each other.
Once it is born, any individual of type $i\in[d]$ gives birth after an exponential time with parameter $\lambda_i>0$ 
to $n_j$ individuals of type $j\in[d]$ with probability $\nu_i(n_1,\dots,n_d)$. Then this individual dies at the same 
time it gives birth. We emphasize that in this model, the probability for the population to become extinct does not 
depend on the rates $\lambda_i$. \\

This model is represented as a plane forest with edge lengths, see Figure \ref{fig2}. 
(In each sibling, we rank individuals of type 1 to the left, then individuals of type 2, and so on.) 
Such a forest will be called a multitype branching forest with edge lengths issued from $x=(x_1,\dots,x_d)$, 
with progeny distribution $\nu:=(\nu_1,\dots,\nu_d)$ and reproduction rates $(\lambda_1,\dots,\lambda_d)$. By 
construction, its discrete time skeleton is a multitype branching (plane) forest, as defined in the previous section, 
with progeny distribution $\nu$, which is independent from the edge lengths. Edge lengths are independent between 
themselves and the length of an edge issued from a vertex of type $i$ follows an exponential distribution with parameter 
$\lambda_i$. We emphasize that the total number of individuals and the total number of mutations in a multitype
branching forest with edge lengths are the same as in its discrete skeleton. Hence, the results of the previous section
can be applied in the present setting.\\

Given a branching forest with edge lengths, as defined above, we denote by $Z=(Z^{(1)},\dots,Z^{(d)})$ the
corresponding multitype branching process, that is for $t\ge0$ and $i\in[d]$, $Z^{(i)}_t$ is the number of 
individuals of type $i$ at time $t$ in the population. (Since no confusion is possible, for the branching process we have 
kept the same notation as in discrete time.) The process $Z$ is a $\mathbb{Z}_+^d$-valued continuous
 time Markov process which satisfies the branching property, i.e., for $\lambda\in\mathbb{R}_+^d$, $t\ge0$ and
$x,y\in\mathbb{Z}_+^d$,
\[\e_{x+y}(e^{-\lambda Z_t})=\e_{x}(e^{-\lambda Z_t})\e_{y}(e^{-\lambda Z_t})\,,\]
where $\p_x$ is the law under which the forest is issued from $x$. In particular, $Z_0=x$, $\p_x$-a.s. 
The process $Z$ actually contains much less information than the original branching forest. In order to preserve 
the essential part of this information, we need to decompose $Z$ as in the following definition. 
\begin{definition}\label{2634}
For $i\neq j$, we denote by $Z^{i,j}_t$ the total number of individuals of type $j$ whose parent has type $i$ and 
who were born before time $t$. For $i=j$, the definition of  $Z^{i,i}_t$ is the same, except that to this number we add the
number of individuals of type $i$ at time $0$ and we subtract the number of individuals of type $i$ who died before 
time $t$. 
\end{definition}
\noindent  The processes $Z^{i,j}$ whose definition should be clear from the example given in Figure \ref{fig2} 
will play a crucial role in our continuous time model. A more formal definition can be found in Section 4.2 of 
\cite{Chaumont}. The interest of these processes is the following straightforward decomposition of the branching 
process $Z=(Z^{(1)},\dots,Z^{(d)})$:
\begin{equation}\label{8254}
Z^{(j)}_t=\sum_{i=1}^dZ^{i,j}_t\,,\;\;\;j\in[d]\,.
\end{equation}
Our model bears on a Lamperti type representation of these processes.  According to Lamperti representation, any one 
dimensional branching process can be expressed as a L\'evy process time changed by some integral functional. In this 
subsection, we will recall from \cite{Chaumont} the extension of this transformation to multitype, continuous time, discrete 
valued branching processes. The latter involves time changed multidimensional compound Poisson processes which we 
now introduce.\\

Since our models of evolution are only concerned with mutations, individuals of type $i$ having 
exactly one child of type $i$ do not present any interest. Hence we can assume without loss of generality that
\[\mbox{$\nu_i(e_i)=0$, for all $i\in[d]$.}\] 
Then let $X=(X^{(1)},\dots,X^{(d)})$, where $X^{(i)}$, $i\in[d]$ are $d$ independent  $\mathbb{Z}^d$-valued 
compound Poisson processes. We assume that $X^{(i)}_0=0$ and that $X^{(i)}$ has rate  $\lambda_i$ and 
jump distribution $\tilde{\nu}_i$ which has been defined in (\ref{4590}). In particular, with the notation 
$X^{(i)}=(X^{i,1},\dots,X^{i,d})$, the process $X^{i,i}$ is a $\mathbb{Z}$-valued, downward skip free, compound 
Poisson process, i.e. $\Delta X_t^{i,i}=X^{i,i}_t-X^{i,i}_{t-}\ge-1$, $t\ge0$, with $X_{0-}=0$ and for all $i\neq j$, the 
process $X^{i,j}$ is an increasing compound Poisson process. We emphasize that in this definition, some of the 
processes $X^{i,j}$, $i,j\in[d]$ can be identically equal to 0.\\ 
    
The following extension of the Lamperti representation to multitype branching processes can be found in
\cite{Chaumont}, see also \cite{cpu} for the case of continuous state multitype branching processes. 

\begin{thm}\label{Lthm} Let us consider a multitype branching forest with edge lengths issued from 
$x=(x_1,\dots,x_d)\in\mathbb{Z}_+^d$, with progeny distribution $\nu:=(\nu_1,\dots,\nu_d)$ and reproduction rates 
$(\lambda_1,\dots,\lambda_d)$. Then the processes $Z^{i,j}$,  $i,j\in[d]$ introduced in Definition $\ref{2634}$ admit the 
following representation:
\begin{equation}\label{3521}
Z^{i,j}_t=\left\{\begin{array}{ll}
X^{i,j}_{\int_0^t Z^{(i)}_s\,ds}\,,\;\;\;t\ge0\,,\;\;\;\mbox{if $i\neq j$,}\\
x_i+X^{i,i}_{\int_0^t Z^{(i)}_s\,ds},\;\;\;t\ge0\,,\;\;\;\mbox{if $i=j$,}
\end{array}\right.
\end{equation}
where the processes, 
\[X^{(i)}=(X^{i,1},X^{i,2},\dots,X^{i,d})\,,\;\;\;i=1,\dots,d\,,\]
are independent $\mathbb{Z}_+^d$ valued compound Poisson processes, with jump distribution 
$(\tilde{\nu}_1,\dots,\tilde{\nu}_d)$ and rates $(\lambda_1,\dots,\lambda_d)$.    
In particular from $(\ref{8254})$ and $(\ref{3521})$, the multitype branching process $Z$ admits the following 
representation, 
\begin{equation}\label{2003}
(Z^{(1)}_t,\dots,Z^{(d)}_t)=x+\left(\sum_{i=1}^dX^{i,1}_{\int_0^t Z^{(i)}_s\,ds},\dots,
\sum_{i=1}^dX^{i,d}_{\int_0^t Z^{(i)}_s\,ds}\right)\,,\;\;\;t\ge0\,.
\end{equation}
\end{thm}

\subsection{Further results on asymptotics of mutations}\label{9523}

\begin{figure}[hbtp]

\hspace*{-0.5cm}{\includegraphics[height=350pt,width=500pt]{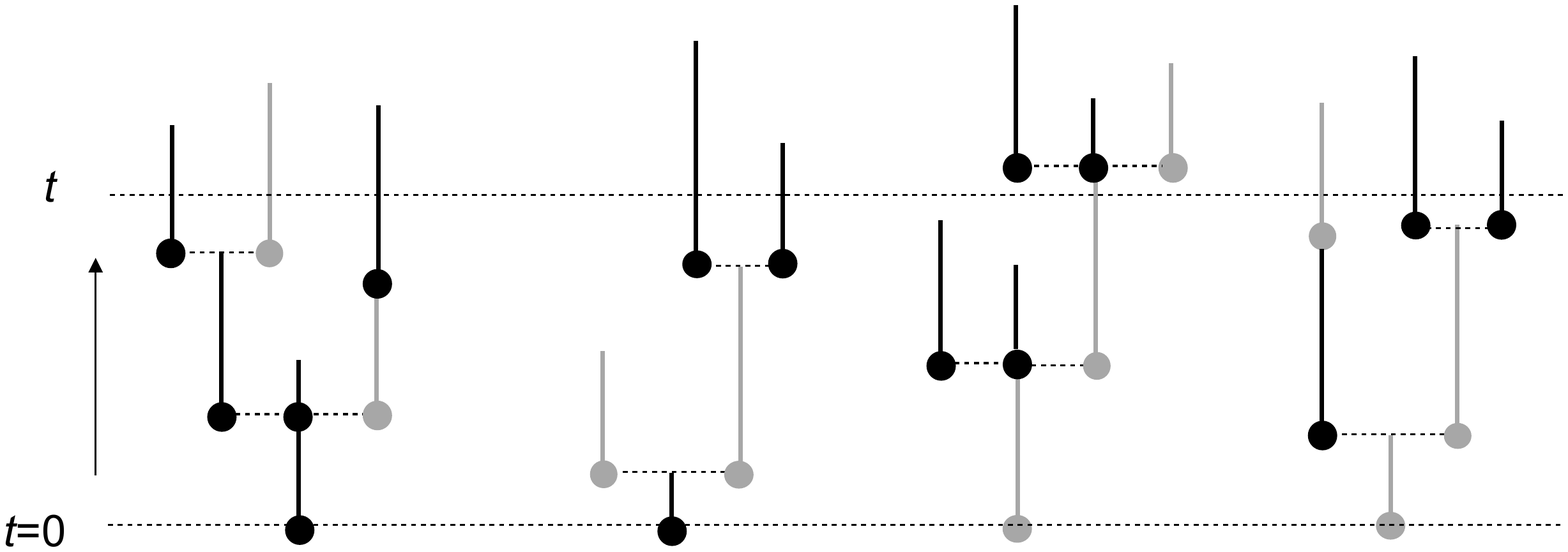}}
 
\vspace{-6.5cm}

\caption{A two type forest with edge lengths issued from $x=(2,2)$. Vertices of type 1 (resp.~2) are represented in black  
(resp.~grey). At time $t$, 
\[Z^{(1)}_t=6,\;\;Z^{(2)}_t=3,\;\; Z^{1,1}_t=-2,\;\;Z^{1,2}_t=5,\;\;Z^{2,1}_t=8,\;\; Z^{2,2}_t=-2,\]
\[\mbox{and}\;\;\;\;M_{1,t}=8,\;\;M_{2,t}=5.\]
}\label{fig2}
\end{figure}

For $i\in[d]$ and $t\ge0$, we will denote by $M_{i,t}$ the total number of mutations of type $i$ which occured
up to time $t$. The definition of this quantity is illustrated on Figure \ref{fig2}. Let us also define a cluster of type $i$ 
as the subtree corresponding to the descendence of type $i$ of an individual of type $i$ which is either a root or 
an individual whose parent as a type different from $i$. Then  $x_i+M_{i,t}$ corresponds to the number 
of clusters of type $i$ in the forest truncated at time $t$.\\

In Proposition \ref{2385}, we describe the asymptotic behaviour of $M_{i,t}$, as $t$ tends to $\infty$ in the case where 
the progeny distribution $\nu$ is primitive and supercritical. To this aim, we will need the joint representation of $M_{i,t}$ 
together with the number $Z_t^{(i)}$ of individuals of type $i$ at time $t$ which is presented in Proposition \ref{Lrep}.

\begin{proposition}\label{Lrep} Recall from Section $\ref{3249}$ the definition of the compound Poisson processes 
$X^{i,j}$, $i,j\in[d]$. Then for any $x=(x_1,\dots,x_d)\in\mathbb{Z}_+^d$, under $\p_x$, the stochastic process 
$\left(Z_t^{(i)},M_{i,t}\right)$ fulfills the following representation,
\[\left(Z_t^{(i)},M_{i,t}\right)=\left(x_i+\sum_{k=1}^d X^{k,i}_{\int_0^t Z^{(k)}_u\,\dif u},
\sum_{k=1, k\neq i}^d X^{k,i}_{\int_0^t Z^{(k)}_u\,\dif u}\right),\quad t\geq 0.\]
\end{proposition}
\begin{proof} This result is a direct consequence of the representation which is recalled in Theorem \ref{Lthm}. Indeed,
recall from Section \ref{3249} the definition of $Z^{i,j}$, then the number of mutations of type $i$ up to time $t$ is 
\[M_{i,t}=\sum_{k\neq i}Z_t^{k,i}\,.\]
The  result follows from identity (\ref{3521}) in Theorem \ref{Lthm}.$\;\;\;\Box$\\   
\end{proof}

Let's us now turn to the limiting behavior of $M_{i,t}$, as $t$ tends to infinity. The next result is concerned with the case 
where $\nu$ is primitive and supercritical. It allows us to evaluate the number of mutations which occured up to time $t$ (or 
equivalently the number of clusters in the forest truncated at time $t$), when $t$ is large.\\ 

Let us define the matrix $A=\Lambda(M-I)$,  where $\Lambda=\mbox{diag}(\lambda_i)$. If $M$ is primitive, then so is $A$ and 
it follows from Perron-Frobenius theory that the eigenvalues $\rho_i$, $i\in[d]$ of $A$ can be arranged so that  
$\rho_1>\mbox{Re}(\rho_2)\ge\dots\ge\mbox{Re}(\rho_d)$. Moreover, $\nu$ is subcritical, critical or supercritical according 
as $\rho_1<0$, $\rho_1=0$ or $\rho_1>0$.  Then a well known result due to \cite{at}, see also Theorem 2, p. 206 in 
\cite{Athreya} asserts that when $\nu$ is non singular and primitive, there exists a nonnegative random variable $W$ such 
that for all $i\in[d]$,
\begin{eqnarray}\label{7392}
\lim_{t\rightarrow\infty}e^{-\rho_1 t}Z_t^{(i)}=v_iW\,,\;\;\mbox{a.s.},
\end{eqnarray}
where $v_i$ is the $i$-th coordinate of the normalized left eigenvector associated with $\rho_1$. 
 
\begin{proposition}\label{2385} Assume that $\nu$ is non singular, primitive and supercritical. Then for all $i\in[d]$,
\[\lim_{t\rightarrow\infty}e^{-\rho_1t}M_{i,t}=K_iW,\ a.s.,\]
where $K_i=v_i(1+(1-m_{ii})(\lambda_i\rho_1)^{-1})$.
\end{proposition}
\begin{proof} We derive from Proposition \ref{Lrep} that,
\[Z_t^{(i)}-M_{i,t}=X^{i,i}_{\int_0^t Z^{(i)}_u\,\dif u}\,,\;\;\;\mbox{a.s.}\]
On the other hand, in the supercritical case, $\rho_1$ is strictly positive. Hence it follows from (\ref{7392}) that
\[\int_0^t Z^{(i)}_u\,\dif u\sim\rho_1^{-1}Wv_ie^{\rho_1 t},\ a.s.,\;\;\ as\ t \to \infty.\]
Then the desired result is a consequence of the latter equivalence and the law of large numbers applied to the 
compound Poisson process $X^{i,i}$. $\;\;\;\Box$\\
\end{proof}

\noindent Under conditions of Propositon \ref{2385}, assume moreover that for some $i\in[d]$, $K_i$ is positive, that is 
\[m_{ii}<1+\lambda_i\rho_1\,,\]
and that for some $j$, $\p_{e_j}(W>0)=1$. Then using Proposition \ref{2385}, we can compare the asymptotic behaviour of the 
number of mutations prior to $t$ with this of ${Z}_t^{(i)}$, under $\p_{e_j}$, that is
\begin{equation}\label{9742}
M_{i,t}\sim K_i{Z}_t^{(i)}\,,\;\;\;\mbox{$\p_{e_j}$-a.s., as $t\rightarrow\infty$}.
\end{equation}
Regarding the condition $\p_{e_j}(W>0)=1$, note that Theorem 2, p. 206 in \cite{Athreya} also asserts that $\p_{e_k}(W>0)>0$, 
for some (hence for all) $k\in[d]$, if and only if 
\begin{equation}\label{9514}
\mbox{$E(\xi_{ij}\log \xi_{ij})<\infty$, for  all $i,j\in[d]$,}
\end{equation}
where $(\xi_{i1},\dots,\xi_{id})$ is a random vector with law $\nu_i$. Moreover, $1-\p_{e_k}(W>0)$ corresponds to
the probability of extinction, when the forest is issued from $e_k$. 

\subsection{Emergence times of mutations} \label{wt}

In this section, we shall assume that mutations are not reversible, that is for all $i=1,\dots,d-1$, individuals of type
$i$ can only have children of type $i$ or $i+1$. In particular $\nu$ is not irreducible. Moreover when giving birth, individuals of 
type $i=1,\dots,d-1$ have at least one child of type $i$ with probability one, and have children of type $i+1$ with positive 
probability. These conditions can be explicited in terms of the progeny distribution $\nu_i$ as follows
\begin{equation}\label{2323}
\left\{\begin{array}{l}
\nu_i(\mathbf{k})>0\;\Rightarrow\;\mbox{$k_j=0$, for $j\notin\{i,i+1\}$,}\\
\sum_{\mathbf{k}\in\mathbb{Z}_+^d:k_{i}=0}\nu_{i}(\mathbf{k})=0\;\;\mbox{and}\;\;
\sum_{\mathbf{k}\in\mathbb{Z}_+^d:k_{i+1}=0}\nu_{i}(\mathbf{k})<1. 
\end{array}\right.
\end{equation}
We are interested in the waiting time until an individual of type $i$ first emerges in the population, that is 
\[\tau_i:=\inf\{t\geq 0:Z^{(i)}_t\geq 1\}\,.\]
The problem of determining a general expression for the law of $\tau_i$ is quite challenging. As far as we know, there is 
no explicit expression for this law in terms of the progeny distribution and the reproduction rates. Various results in this 
direction can be found in \cite{Serra}, \cite{Serra-Haccou},  \cite{Durrett2010} and \cite{Alexander} for instance. Most of 
them provide approximations of this law, using martingale convergence theorems \cite{Durrett2010} or through 
numerical methods for the inversion of the generating function \cite{Alexander}. In Proposition \ref{4171} we first give 
a relationship between the successive emergence times $\tau_2,\tau_3,\dots$ in terms of the underlying compound 
Poisson process in the Lamperti representation of $Z$. We also characterize the joint law under $\p_{e_{i-1}}$ of the time 
$\tau_i$ and the number of individuals of type $i-1$ at this time. In Theorem \ref{9241}  we derive an approximation of  the 
time $\tau_i$, under $\p_{e_1}$, as the mutation rate of type $k$ increases faster than that of type $k-1$, for all $k=3,\dots,i$. 
Then in Corollary \ref{2482} we focus on a case where these law can be explicited.\\

In the following developments, we use the notation of Section \ref{3249} from which we recall the Lamperti representation 
of the multitype branching process $Z=(Z^{(1)},\dots,Z^{(d)})$ in terms of the compound Poisson processes $X^{(i)}$. 
Let us also introduce a few more notation. For $i,j\in[d]$, we denote by $\lambda_{i,j}$ the parameter of the
compound Poisson process $X^{i,j}$, that is
\[\lambda_{i,j}:=\mbox{$\lambda_i\left(1-\sum_{\mathbf{k}\in\mathbb{Z}_+^d:k_j=0}\tilde{\nu}_{i}(\mathbf{k})
\right).$}\]
Note that from our assumptions (\ref{2323}), for all $i=1,\dots,d-1$, $\lambda_{i,i+1}>0$ and for $j\notin\{i,i+1\}$, 
$\lambda_{i,j}=0$, that is $X^{i,j}$ is identically equal to 0. In particular, $\lambda_i=\lambda_{i,i}+\lambda_{i,i+1}$, for
$i\le d-1$ and $\lambda_d=\lambda_{d,d}$.
The parameter $\lambda_{i,i+1}$ will be call the mutation rate of type $i+1$. For $i\ge2$, let 
\[\gamma_i:=\inf\{t:X^{i-1,i}_t\ge1\}\]
be the time of the first jump by the process $X^{i-1,i}$ and note that this time is exponentially distributed with parameter
$\lambda_{i-1,i}$. 

\begin{proposition}\label{4171}   
Assume that $(\ref{2323})$ holds and define $Z^{0,1}$ as the process identically equal to $1$ and set $\tau_1=0$.
\begin{itemize}
\item[$1.$]  For $i=2,\dots,d$, the emergence time $\tau_i$ of type $i$ admits the following representation under 
$\p_{e_1}$,
\begin{equation}\label{3975}
\tau_i=\tau_{i-1}+\int_{0}^{\gamma_i}
\frac{1}{X^{i-1,i-1}_s+Z^{i-2,i-1}_{\kappa_{i-1}(s)}}\,\dif s,\;\;\;\mbox{$\p_{e_1}$-a.s.,}
\end{equation}
where $\kappa_{i-1}$ is the right continuous inverse of the functional $t\mapsto\int_{0}^{t}Z^{(i-1)}_s\,\dif s\), 
i.e. $\kappa_{i-1}(t)=\inf\{s>0:\int_{0}^{s}Z^{(i-1)}_u\,\dif u>t\}$.
\item[$2.$] Under $\p_{e_{i-1}}$, the joint law of the emergence time $\tau_i$ 
of type $i$ together with the number of individuals of type $i-1$ in the population at time $\tau_i$ admits the following 
representation,
\begin{equation}\label{3976}
(\tau_{i},Z^{(i-1)}_{\tau_i})\overset{\textup{\textrm{(d)}}}{=}\left(\int_{0}^{\gamma_i}\frac{\dif s}{1+
X_s^{i-1,i-1}},1+X^{i-1,i-1}_{\gamma_i}\right).
\end{equation}
\item[$3.$]  Let us define $\theta_k=\int_{0}^{\gamma_k}\frac{1}{X^{k-1,k-1}_{s}+1}\,\dif s$,
for $k\ge2$.Then the random variables $\theta_k$, $k\ge2$ are independent and for $i=2,\dots,d$,
\begin{equation}\label{1086}
\p_{e_1}(\tau_i>t)\le P\left(\sum_{k=2}^{i}\theta_k>t\right)\,,\;\;\mbox{for all $t>0$\,.}
\end{equation}
\end{itemize}
\end{proposition}
\begin{proof} Since $X^{i,j}$ is identically equal to 0 whenever $j\notin\{i,i+1\}$, then under $\p_{e_1}$, the 
representation (\ref{2003})  admits the simpler form
\begin{equation}\label{1531}
(Z^{(1)}_t,\dots,Z^{(d)}_t)=e_1+\left(X^{1,1}_{\int_0^t Z^{(1)}_s\,ds},X^{2,2}_{\int_0^t Z^{(2)}_s\,ds}+
X^{1,2}_{\int_0^t Z^{(1)}_s\,ds},\dots,X^{d,d}_{\int_0^t Z^{(d)}_s\,ds}+X^{d-1,d}_{\int_0^t Z^{(d-1)}_s\,ds}\right).
\end{equation}
In particular, for $i=2,\dots,d$,
\[Z^{(i)}_t=X^{i,i}_{\int_0^t Z^{(i)}_s\,ds}+X^{i-1,i}_{\int_0^t Z^{(i-1)}_s\,ds},\;\;t\ge0.\]
Since $X^{i,i}_0=0$, for $i\ge2$, we see that the time $\tau_i$ corresponds to the first hitting time of level 1 by the 
process  $t\mapsto X^{i-1,i}_{\int_0^t Z^{(i-1)}_s\,ds}$, that is 
\begin{equation}\label{3589}
\tau_i=\kappa_{i-1}(\gamma_i),
\end{equation}
where $\gamma_i$ has been defined as the time of the first jump of the process $X^{i-1,i}$. For $t$ such that
$\kappa_{i-1}(t)<\infty$, we have  $t=\int_0^{\kappa_{i-1}(t)}Z^{(i-1)}_s\,\dif s$, so that 
$\dif t=Z^{(i-1)}_{\kappa_{i-1}(t)}\dif \kappa_{i-1}(t)$, and since  $\kappa_{i-1}(0)=\tau_{i-1}$, we obtain
\begin{eqnarray}
\kappa_{i-1}(t)&=&\tau_{i-1}+\int_0^t\frac{\dif s}{Z^{(i-1)}_{\kappa_{i-1}(s)}}\nonumber\\
&=&\tau_{i-1}+\int_0^t\frac{\dif s}{X^{i-1,i-1}_s+X^{i-2,i-1}_{\int_0^{\kappa_{i-1}(s)} Z^{(i-2)}_u\,\dif u}}\,.\label{4278}
\end{eqnarray}
The latter identity together with (\ref{3589}) prove identity (\ref{3975}).

The second part of the proposition is easily derived from the same arguments. More specifically, it follows from 
(\ref{3589}) and  the following identities
\[Z^{(i-1)}_t=1+X^{i-1,i-1}_{\int_0^t Z^{(i-1)}_s\,ds}\;\;\mbox{and}\;\;\kappa_{i-1}(t)=\int_{0}^{t}\frac{\dif s}{1+
X_s^{i-1,i-1}}\,,   \;\;t\ge0\,,\]
which hold $\p_{e_{i-1}}$-a.s.

Independence between the variables $\theta_k$, $k\ge2$ is a direct consequence of the independence between the 
processes $X^{(i)}$, $i\in[d]$. We derive from the representation of $\tau_i$ in part 1. of this proposition that
\begin{equation}\label{2431}
\tau_i=\sum_{k=2}^i \int_{0}^{\gamma_k}\frac{1}{X^{k-1,k-1}_s+
X^{k-2,k-1}_{\int_0^{\kappa_{k-1}(s)} Z^{(k-2)}_u\,\dif u}}\,\dif s,\;\;\;\mbox{a.s.}
\end{equation}
Note that since $\kappa_{k-1}(0)=\tau_{k-1}$, then from (\ref{3589}), for all $k\ge2$, 
$\int_0^{\kappa_{k-1}(0)} Z^{(k-2)}_u\,\dif u=\gamma_{k-1}$, so that by definition of $\gamma_{k-1}$,
\begin{equation}\label{2381}
X^{k-2,k-1}_{\int_0^{\kappa_{k-1}(0)} Z^{(k-2)}_u\,\dif u}=X^{k-2,k-1}_{\gamma_{k-1}}\ge1,\;\;\;\mbox{a.s.}
\end{equation}
Besides, since $s\mapsto X^{k-2,k-1}_{\int_0^{\kappa_{k-1}(s)} Z^{(k-2)}_u\,\dif u}$ are increasing processes, then
inequality (\ref{1086}) is a direct consequence of identities (\ref{2431}) and (\ref{2381}).
$\;\;\;\Box$\\
\end{proof} 
\noindent Note that the law of $\theta_k$ or equivalently, the law of $\tau_k$ under $\p_{e_{k-1}}$ can be made 
explicit in some instances through its Laplace transform, see Corollary \ref{2482} below.\\

For the remainder of this section we will assume moreover that at each mutation, individuals of type $i$ 
do not give birth to more than one child of type $i+1$ in a same litter. More specifically, assumptions (\ref{2323}) 
are replaced by,
\begin{equation}\label{4872}
\left\{\begin{array}{l}
\nu_i(\mathbf{k})>0\;\Rightarrow\;\mbox{$k_{i+1}=0$ or 1 and $k_j=0$, for $j\notin\{i,i+1\}$,}\\
\sum_{\mathbf{k}\in\mathbb{Z}_+^d:k_{i}=0}\nu_{i}(\mathbf{k})=0\;\;\mbox{and}\;\;
\sum_{\mathbf{k}\in\mathbb{Z}_+^d:k_{i+1}=0}\nu_{i}(\mathbf{k})<1\,.
\end{array}\right.
\end{equation}
In particular, under these assumptions, the process $X^{i,i+1}$ is a standard Poisson process. 
Then we  will need the next lemma in order to derive our main result on the estimation of the time $\tau_i$,
as the mutation rates $\lambda_{k-1,k}$, $k=2,\dots,d$ grow faster. 

\begin{lemma}\label{4626} Assume that $(\ref{4872})$ holds, let $k\ge3$ and fix $\lambda_{1,2}>0$, then
\begin{eqnarray}\label{5435}
&&\p_{e_1}\big(X^{k-2,k-1}_{\int_0^{\kappa_{k-1}(\gamma_k)} Z^{(k-2)}_u\,\dif u}=1\big)\longrightarrow1\,,\nonumber\\
&&\qquad\qquad\qquad\qquad\qquad\mbox{as $\lambda_{n-2,n-1}/\lambda_{n-1,n}\rightarrow0$, for $n=3,\dots,k$.}
\end{eqnarray}
\end{lemma}

\begin{proof} First set $\gamma_{k-1}^{(1)}=\inf\{t>\gamma_{k-1}:X^{k-2,k-1}_t=2\}$ and note that
\begin{eqnarray*}
\mbox{$\{X^{k-2,k-1}_{\int_0^{\kappa_{k-1}(\gamma_k)} Z^{(k-2)}_u\,\dif u}=1\}$}&=&
\mbox{$\{\int_0^{\kappa_{k-1}(\gamma_k)} Z^{(k-2)}_u\,\dif u<\gamma_{k-1}^{(1)}\}$}\\
&=&\mbox{$\{\kappa_{k-1}(\gamma_k)<\kappa_{k-2}(\gamma_{k-1}^{(1)})\}.$}
\end{eqnarray*}
It is easy to check that $\kappa_{k-2}(\gamma_{k-1}^{(1)})=\tau_{k-1}^{(1)}$, where
\[\tau_{k-1}^{(1)}:=\inf\{t>\tau_{k-1}:Z^{k-2,k-1}_t-Z^{k-2,k-1}_{\tau_{k-1}}=1\}\,.\] 
(Note that from our assumptions $Z^{k-2,k-1}_{\tau_{k-1}}=1$ and $Z^{k-2,k-1}_{\tau_{k-1}^{(1)}}=2$, $\p_{e_1}$-a.s.) 
So from (\ref{3589}), we have showed that 
\begin{equation}\label{1237}
\mbox{$\{X^{k-2,k-1}_{\int_0^{\kappa_{k-1}(\gamma_k)} Z^{(k-2)}_u\,\dif u}=1\}$}=\{\tau_{k}<\tau_{k-1}^{(1)}\}\,.
\end{equation}
The event $\{\tau_{k}<\tau_{k-1}^{(1)}\}$ means that before the first time when an individual of type $k$ appears in the 
population, there has been only one birth of type $k-1$. From the Markov property applied at time $\tau_{k-1}$, we have 
\begin{equation}\label{5583}
\p_{e_1}(\tau_{k}\le \tau_{k-1}^{(1)})=\int\p_{z}(\tau_{k}\le \tau_{k-1}^{(1)})\p_{e_1}(Z_{\tau_{k-1}}\in dz)\,.
\end{equation}
The support in the integral of (\ref{5583}) is included in the set $\{z:z_{k-1}=1\}$, so
from (\ref{1237}), (\ref{5583}) and the Lebesgue theorem of dominated convergence, all we need to prove is 
\begin{equation}\label{3248}
\p_{z}(\tau_{k}\le \tau_{k-1}^{(1)})\rightarrow1,\;\;\;\mbox{as 
$\lambda_{n-2,n-1}/\lambda_{n-1,n}\rightarrow0$, for $n=3,\dots,k$,}
\end{equation}
for all $z$ such that $z_{k-1}=1$. (Note that if $z$ is such that $z_1=\dots=z_{k-2}=0$, or such that $z_k\ge1$, then it is clear 
that $\p_{z}(\tau_{k}\le \tau_{k-1}^{(1)})=1$, since in the first case $Z^{k-2,k-1}$ is identically equal to 0, so that 
$\tau_{k-1}^{(1)}=\infty$, $\p_z$-a.s. and in the second case, $\tau_k=0$, $\p_z$-a.s.)  

Let $z$ be such that $z_{k-1}=1$. Without loss of generality we can assume that $z_i\ge1$, for $i=1,\dots k-2$.
For $i=1,\dots k-1$, let us denote by $U_i$ the first time that the lineage of one of the 
$z_{k-i}$ initial individuals of type $k-i$ gives birth to an individual of type $k-i+1$. Then from the branching property, under 
$\p_z$, the r.v.'s $U_i$ are independent and from part 2. of Proposition \ref{4171}, $U_i$ has the same law as 
$\int_0^{\gamma_{k-i+1}}\frac{\dif s}{X_s^{k-i,k-i}+z_{k-i}}$. Then set $Y_s^{(i)}:=X_s^{k-i,k-i}+z_{k-i}$ and note the inclusions, 
\[\left\{\gamma_k\le \min\left(\gamma_{k-1}/Y^{(2)}_{\gamma_{k-1}},\dots,
\gamma_{2}/Y^{(k-1)}_{\gamma_2}\right)\right\}\subset\{U_1\le \min(U_2,\dots,U_{k-1})\}\subset\{\tau_{k}\le \tau_{k-1}^{(1)}\}\,,\]
which imply the inequality,
\[P\big(\gamma_k/\gamma_{k-1}\le \min\big(1/Y^{(2)}_{\gamma_{k-1}},\gamma_{k-2}/(\gamma_{k-1}Y^{(3)}_{\gamma_{k-2}}),
\dots,\gamma_{2}/(\gamma_{k-1}Y^{(k-1)}_{\gamma_2})\big)\big)\le \p_z(\tau_{k}\le \tau_{k-1}^{(1)})\,.\]
But when $\lambda_{n-2,n-1}/\lambda_{n-1,n}\rightarrow0$, for $n=3,\dots,k$, the parameter $\lambda_{1,2}>0$ being fixed, 
we necessarily have $\lim\lambda_{n-1,n}=\infty$, for $n=3,\dots,k$. Hence $\gamma_k/\gamma_{k-1}$ converges in 
probability toward 0, $1/Y^{(2)}_{\gamma_{k-1}}$ converges in probability toward $1/z_{k-2}$ and 
$\gamma_{n-1}/(\gamma_{n}Y^{(k-n+2)}_{\gamma_{n-1}})$, for $n=3,\dots,k-1$ converge in probability toward $+\infty$. 
Therefore, the left hand side of the above inequality tends to 1, which proves (\ref{3248}) and the lemma is proved. 
$\;\;\;\Box$\\
\end{proof} 

\noindent In the following theorem, the assumption $\frac{\lambda_{k-1,k}}{\lambda_{k,k+1}}\rightarrow0$ is quite adapted to 
several biological models such as  cancer growth, for instance. Indeed, cancer is often the result of a series of successive 
mutations, \cite{Iwasa2}, \cite{Durrett2010}, \cite{Durrett}. Each new mutation is itself more unstable than the previous ones, 
and in particular, 
the successive mutation rates can increase very fast. It would interesting to study the asymptotic behavior of $\tau_i$, when 
$\frac{\lambda_{k,k}}{\lambda_{k+1,k+1}}\rightarrow0$, that is when the intrinsic reproduction rates increase very fast. 
This assumption also fits to the model of cancer, since mutations are always more sensitive to proliferate.

\begin{thm}\label{9241} Assume that $(\ref{4872})$ holds. Recall the definition of $\theta_k$ in Proposition $\ref{4171}$ and  
let us fix $\lambda_{1,2}>0$, then under $\p_{e_1}$,
\[\frac{\tau_i}{\sum_{k=2}^{i}\theta_k}
\stackrel{P}{\longrightarrow}1\,,\;\;\mbox{as $\frac{\lambda_{k-2,k-1}}{\lambda_{k-1,k}}\rightarrow0$, for $k=3,\dots,i$.}\]
Besides, the expectation of $\tau_i$ fulfills the following approximation:
\[\e_{e_1}(\tau_i)\sim\sum_{k=2}^iE(\theta_k)\,,\;\;
\mbox{as $\frac{\lambda_{k-2,k-1}}{\lambda_{k-1,k}}\rightarrow0$, for $k=3,\dots,i$.}\]
\end{thm}
\begin{proof} Since $s\mapsto X^{k-2,k-1}_{\int_0^{\kappa_{k-1}(s)} Z^{(k-2)}_u\,\dif u}$ are increasing processes, then 
from (\ref{2381}), $\p_{e_1}$-almost surely on the set $\{X^{k-2,k-1}_{\int_0^{\kappa_{k-1}(\gamma_k)} Z^{(k-2)}_u\,\dif u}=1\}$,
we have
\[\int_{0}^{\gamma_k}\frac{1}{X^{k-1,k-1}_s+X^{k-2,k-1}_{\int_0^{\kappa_{k-1}(s)} Z^{(k-2)}_u\,\dif u}}\,\dif s=
\int_{0}^{\gamma_k}\frac{1}{X^{k-1,k-1}_s+1}\,\dif s\,.\]
Hence it follows from Lemma \ref{4626} that for fixed $\lambda_{1,2}>0$, as $\lambda_{n-2,n-1}/\lambda_{n-1,n}\rightarrow0$, 
for all $n=3,\dots,k$, 
\[\left(\int_{0}^{\gamma_k}\frac{1}{X^{k-1,k-1}_s+1}\,\dif s\right)^{-1}\int_{0}^{\gamma_k}\frac{1}{X^{k-1,k-1}_s
+X^{k-2,k-1}_{\int_0^{\kappa_{k-1}(s)} Z^{(k-2)}_u\,\dif u}}\,\dif s\stackrel{P}{\longrightarrow}1\,,\]
and the first part of the theorem is easily derived from this convergence and (\ref{3975}) (or equivalently (\ref{2431})). 

In order to prove the second part, let us first set 
\[H_k:=\int_{0}^{\gamma_k}\frac{1}{X^{k-1,k-1}_s+X^{k-2,k-1}_{\int_0^{\kappa_{k-1}(s)} Z^{(k-2)}_u\,\dif u}}\,\dif s\;\;\mbox{and}
\;\;A_k:=\{X^{k-2,k-1}_{\int_0^{\kappa_{k-1}(\gamma_k)} Z^{(k-2)}_u\,\dif u}=1\}\,.\]
Then from (\ref{3975}), $\e_{e_1}(\tau_i)=\sum_{k=2}^i\e_{e_1}(H_k)$, so it suffices to prove that for all $k=2,\dots,i$, 
\begin{equation}\label{2363}
\e_{e_1}(H_k)\sim E(\theta_k),\;\;\mbox{as $\lambda_{n-2,n-1}/\lambda_{n-1,n}\rightarrow0$, for $n=3,\dots,k$.}
\end{equation}
Observe that $\e_{e_1}(H_k)=E(\theta_k{\bf 1}_{A_k})+\e_{e_1}(H_k{\bf 1}_{A_k^c})$. Moreover, $\e_{e_1}(H_k{\bf 1}_{A_k^c})
\le E(\theta_k{\bf 1}_{A_k^c})$. Then to obtain (\ref{2363}), it is enough to prove that 
\begin{equation}\label{8441}
\frac{E(\theta_k{\bf 1}_{A_k^c})}{E(\theta_k)}\rightarrow0,\;\;\mbox{as $\lambda_{n-2,n-1}/\lambda_{n-1,n}\rightarrow0$, for 
$n=3,\dots,k$.}
\end{equation}
But for any $p,q\ge1$, such that $p^{-1}+q^{-1}=1$, we have from Holder inequality 
$E(\theta_k{\bf 1}_{A_k^c})\le E(\theta_k^p)^{1/p}P(A_k^c)^{1/q}$. Moreover, we clearly have 
$E(\theta_k^p)^{1/p}\sim 1/\lambda_{k-1,k}$, as $\lambda_{k-1,k}\rightarrow\infty$. Hence, (\ref{8441}) is satisfied thanks to 
Lemma \ref{4626}. $\;\;\;\Box$\\
\end{proof} 

We end this section with an example where the distribution of $\tau_i$ can be estimated a bit more specifically. 
We consider the case of binary fission with mutations, where each individual of type $i$ can give birth to either two individuals 
of type $i$ or one individual of type $i$ and one individual of type $i+1$. In particular, all jumps of $Z^{i,i}$ have size $1$ and 
$X^{i,i}$ is a Poisson process with parameter $\lambda_{i,i}$. 

\begin{corollary}\label{2482} With the above assumtions, the law of $\tau_i$ can be specified as follows.
\begin{itemize}
\item[$1.$] Under $\p_{e_{i-1}}$, the Laplace transform of $\tau_i$ is expressed as,
\[\e_{e_{i-1}}(e^ {-\alpha \tau_i})=
\lambda_{i-1,i}\sum_{n\geq 0}\frac{\lambda_{i-1,i-1}^n}{\prod_{k=0}^n(\alpha_k+\dots+\alpha_{n}+\bar{\alpha}_{n+1})},
\;\;\;\alpha\ge0\,,\]
where $\alpha_0=0$, $\alpha_k=\frac{\alpha}{k(k+1)}$ and  $\bar{\alpha}_{k}=\lambda_{i-1}+\frac\alpha{k}$, for $k\ge1$.
\item[$2.$] The expectation of $\tau_i$ is given by 
$\e_{e_{i-1}}(\tau_i)=\frac{1}{\lambda_{i-1,i}\lambda_{i-1,i-1}}\ln\frac{\lambda_{i-1}}{\lambda_{i-1,i}}$.
In particular, for fixed $\lambda_{1,2}>0$, under $\p_{e_1}$, the expectation of $\tau_i$ fulfills the following approximation:
\[\e_{e_1}(\tau_i)\sim\sum_{k=2}^i\lambda_{k-1,k}^{-2}\,,\;\;\;
\mbox{as $\frac{\lambda_{k-2,k-1}}{\lambda_{k-1,k}}\rightarrow0$, for $k=3,\dots,i$.}\]
\end{itemize}
\end{corollary}

\begin{proof} From part 2. of Proposition \ref{4171}, for all $\beta\ge0$,
\begin{eqnarray}
\e_{e_{i-1}}(e^{-\alpha\tau_i})&=&E\left(e^{-\alpha\int_{0}^{\gamma_i}\frac{1}{1+X_s^{i-1,i-1}}\dif s}\right)\nonumber\\
&=&\lambda_{i-1,i}\int_{0}^{+\infty}E\left(e^{-\alpha\int_{0}^{x}\frac{1}{1+X_s^{i-1,i-1}}\dif s}\right)
e^{-\lambda_{i-1,i}x}\dif x\,.\label{6423}
\end{eqnarray}
Under $\p_{e_{i-1}}$, $X^{i-1,i-1}$ is a standard Poisson process with parameter $\lambda_{i-1,i-1}$ starting 
at $0$. So if we denote by $(J_n)_{n\ge1}$ the sequence of jump times of $X^{i-1,i-1}$ and set $J_0=0$, then developing 
the expression $E\left(e^{-\alpha\int_{0}^{x}\frac{1}{1+X_s^{i-1,i-1}}\dif\, s}\right)$, we obtain with the convention that 
$\sum_{k=0}^{-1}=0$,
\begin{align*}
E\left(e^{-\alpha\int_{0}^{x}\frac{1}{1+X_s^{i-1,i-1}}\dif\, s}\right)
&=\sum_{n\geq 0}E\left(X_x^{i-1,i-1}=n, e^{-\alpha\left(\frac{x-J_n}{n+1}+\sum_{k=0}^{n-1}\frac{J_{k+1}-J_k}{k+1}\right)}\right)\\
&=e^{-(\alpha+\lambda_{i-1,i-1})x}+\sum_{n\geq 1}e^{-\lambda_{i-1,i-1} x}\frac{(\lambda_{i-1,i-1}x)^n}{n!}\\
&\qquad\qquad\times\int_{0\leq x_1\leq \dots \leq x_n\leq x}
e^{-\alpha\left(\frac{x}{n+1}+\sum_{k=1}^{n}\frac{x_k}{k(k+1)}\right)}\frac{n!}{x^n}\dif x_1\dots\dif x_n\\
&=e^{-(\alpha+\lambda_{i-1,i-1})x}+\sum_{n\geq 1}\lambda_{i-1,i-1}^ne^{-\left(\lambda_{i-1,i-1}+\frac{\alpha}{n+1}\right)x}\\
&\qquad\qquad\times\int_{0\leq x_1\leq \dots \leq x_n\leq x}e^{-\alpha\sum_{k=1}^{n}\frac{x_k}{k(k+1)}}\dif x_1\dots\dif x_n.
\end{align*}
Then coming back to expression (\ref{6423}), we obtain with the convention that $\sum_{k=1}^{0}=0$,
\[\e_{e_{i-1}}(e^{-\alpha\tau_i})=\lambda_{i-1,i}\sum_{n\geq 0}\lambda_{i-1,i-1}^n
\int_{0\leq x_1\leq \dots \leq x_{n+1}}e^{-(\bar{\alpha}_{n+1}x_{n+1}+\sum_{k=1}^{n}\alpha_kx_k)}\dif x_1\dots\dif x_{n+1},\]
where $\alpha_1,\dots,\alpha_{n},\bar{\alpha}_{n+1}$ are defined in the satement. (Here we used the 
fact that $\lambda_{i-1}=\lambda_{i-1,i}+\lambda_{i-1,i-1}$.) The computation of the integral is easily done.

Then using again part 2. of Proposition \ref{4171}, we obtain the expectation of $\tau_i$ under $\p_{e_{i-1}}$, after easy 
computations,
\begin{align*}
\e_{e_{i-1}}(\tau_i)
&=\int_0^{+\infty} \dif x \lambda_{i-1,i}e^{-\lambda_{i-1,i}x}\int_0^x e^{-\lambda_{i-1,i-1}s}\sum_{k\geq 0}
\frac{(\lambda_{i-1,i-1}s)^k}{(k+1)!} \dif s\\
&=\frac{1}{\lambda_{i-1,i}\lambda_{i-1,i-1}}\ln\frac{\lambda_{i-1}}{\lambda_{i-1,i}}\,.\\
\end{align*}
We conclude from Theorem \ref{9241}.
$\;\;\;\Box$\\
\end{proof}

\end{document}